\documentclass[a4paper,10pt]{article}
\usepackage[french]{babel}
\usepackage[utf8]{inputenc}
\usepackage{setspace}
\usepackage{cite}
\usepackage{color}
\usepackage{amssymb}
\usepackage{amsmath}
\usepackage{amsthm}
\usepackage{graphicx}
\usepackage{enumerate}
\usepackage{eufrak}

\def\CA{{\mathcal A}}

\def\CH{{\mathcal H}}

\def\CP{{\mathcal P}}

\def\pa{{\partial}}

\newcommand{\sca}[2]{\langle#1,#2\rangle}

\newcommand{\beq}{\begin{equation}}
\newcommand{\eeq}{\end{equation}}

\renewcommand{\Re}{{\rm Re\,}}
\renewcommand{\Im}{{\rm Im\,}}
\newtheorem{theorem}{Theorem}[section]

\newtheorem{corollary}[theorem]{Corollary}

\makeatletter
\@addtoreset{equation}{section}

\makeatother

\title{Spectral instability of some non-selfadjoint anharmonic oscillators\\
Instabilité spectrale de certains oscillateurs anharmoniques non-autoadjoints}

\author{Raphaël HENRY}

\begin{document}
\maketitle

\begin{abstract}
Notre objectif est de mettre en évidence l'instabilité spectrale de certains opérateurs différentiels non-autoadjoints, via
l'étude de la croissance des normes des projecteurs spectraux $\Pi_n$ associés à leurs valeurs propres.
Nous nous intéressons à certains \textit{oscillateurs anharmoniques} 
\[\label{oscanh}
\CA(m,\theta) = -\frac{d^2}{dx^2}+e^{i\theta}|x|^m
\]
avec
\[\label{argz}
|\theta|<\min\left\{\frac{(m+2)\pi}{4},\frac{(m+2)\pi}{2m}\right\},
\]
définis sur $L^2(\mathbb{R})$. Nous 
étendons les résultats de \cite{Dav2} et \cite{DavKui} en donnant
 un développement asymptotique de la norme des projecteurs spectraux associés aux
grandes valeurs propres pour les opérateurs $\CA(1,\theta)$ et $\CA(2k,\theta)$, $k\geq1$.\\

The purpose of this Note is to highlight the spectral instability of some non-selfadjoint differential operators, by 
studying the growth rate of the norms of the spectral projections $\Pi_n$ associated with their eigenvalues.
More precisely, we are concerned with some \textit{anharmonic oscillators}
\[\label{oscanh}
\CA(m,\theta) = -\frac{d^2}{dx^2}+e^{i\theta}|x|^m
\]
with
\[\label{argz}
|\theta|<\min\left\{\frac{(m+2)\pi}{4},\frac{(m+2)\pi}{2m}\right\},
\]
defined on $L^2(\mathbb{R})$. We get asymptotic expansions for the norm of the spectral projections associated with
the large eigenvalues
of $\CA(1,\theta)$ and $\CA(2k,\theta)$, $k\geq1$, extending the results of \cite{Dav2} and \cite{DavKui}.
\end{abstract}

\section{Spectral instability and pseudospectra}
It is well-known that the spectral theorem implies some control of stability for the spectrum of selfadjoint operators : if $\CA$
is a selfadjoint operator acting on the Hilbert space $\CH$, the spectrum of its perturbations
$\CA+\varepsilon \mathcal{B}$, with $\varepsilon>0$ and any $\mathcal{B}\in\mathcal{L}(\CH)$, 
$\|\mathcal{B}\|\leq1$, lies entirely inside an 
$\varepsilon$-neighborhood of the spectrum $\sigma(\CA)$.
In other words, the norm of
the resolvent of $\CA$ near the spectrum blows up like the inverse distance to the spectrum. It has also been known
for several years (see for instance \cite{TrEm}) that such a behavior could not be expected in general in the case of 
non-selfadjoint operators.
One can understand it thanks to the notion of \textit{$\varepsilon$-pseudospectra} of an operator $\CA$, 
defined as the family of sets $\sigma_\varepsilon(\CA)$, indexed by $\varepsilon>0$,
\[
 \sigma_\varepsilon(\CA)  =  \left\{\xi\in\rho(\CA) : \|(\CA-\xi)^{-1}\|>\frac{1}{\varepsilon}\right\}\cup\sigma(\CA).
\]
The link between spectral instability and pseudospectra appears more clearly in the following equivalent formulation, which
is a weak version of the Roch and Silbermann theorem \cite{RoSi} :
\[
 \sigma_\varepsilon(\CA) =  \bigcup_{
\tiny{\begin{array}{c}\omega\in\mathcal{L}(\CH),\\ \|\mathcal{B}\|\leq1
\end{array}}
} \sigma(\CA+\varepsilon\mathcal{B})
\]
(See also \cite{Sjo} and references therein).\\
In the following, we deal with the instability indices associated with an isolated eigenvalue
$\displaystyle{\lambda\in\sigma(\CA)}$. The instability index associated with $\lambda$ is defined as
\[
 \kappa(\lambda) = \|\Pi(\lambda)\|,
\]
where $\Pi(\lambda)$ denotes the spectral projection associated with $\lambda$. Of course $\kappa(\lambda)\geq1$ in any case,
and $\kappa(\lambda) = 1$ when $\CA$ is selfadjoint.
These numbers $\kappa(\lambda)$ are closely related to the size of $\varepsilon$-pseudospectra around $\lambda$. 
Indeed, if $\sigma_\varepsilon^\lambda$ denotes the connected component of $\sigma_\varepsilon(\CA)$ containing $\lambda$,
 and if we assume for simplicity that
$\sigma_\varepsilon^\lambda\cap\sigma(\CA) = \{\lambda\}$ and $\sigma_\varepsilon^\lambda$ is bounded,
then the perimeter $|\pa\sigma_\varepsilon^\lambda|$ of $\sigma_\varepsilon^\lambda$ satisfies (see \cite{AslDav})
\beq
|\pa\sigma_\varepsilon^\lambda|\geq2\pi\varepsilon\kappa(\lambda).
\eeq
In the finite dimensional setting at least, instability indices give a better description of pseudospectra
: if $\CA\in\mathcal{M}_n(\mathbb{C})$ is a diagonalisable matrix with distinct eigenvalues 
$\lambda_1,\dots,\lambda_n$, Embree and Trefethen
show \cite{TrEm} that
there exists $\varepsilon_0>0$ such that, for all $\varepsilon<\varepsilon_0$,
\beq\label{incluPseudo}
 \bigcup_{\tiny{\lambda_k\in\sigma(\CA)}}D(\lambda_k,\varepsilon\kappa(\lambda_k)+\mathcal{O}(\varepsilon^2))\subset
\sigma_\varepsilon(\CA)\subset
\bigcup_{\tiny{\lambda_k\in\sigma(\CA)}}D(\lambda_k,\varepsilon\kappa(\lambda_k)+\mathcal{O}(\varepsilon^2)).
\eeq
In the case of an infinite dimensional space, the validity of this statement should be investigated, as well as the
dependance on $\lambda_k$ of the $\mathcal{O}(\varepsilon^2)$ terms.\\
In the following, we study the instability indices of simple non-selfadjoint differential operators introduced by Davies in
\cite{Dav2}, for which the instability phenomenon described above will appear clearly.
Let us define the \textit{anharmonic oscillators} 
\beq\label{oscanh}
\CA(m,\theta) = -\frac{d^2}{dx^2}+e^{i\theta}|x|^m
\eeq
with
\beq\label{argz}
|\theta|<\min\left\{\frac{(m+2)\pi}{4},\frac{(m+2)\pi}{2m}\right\},
\eeq
defined on $L^2(\mathbb{R})$ in \cite{Dav2} by taking the closure of the associated quadratic form defined on
$\mathcal{C}_0^\infty(\mathbb{R})$,
which is sectorial if $\theta$ satisfies (\ref{argz}).
According to \cite{Dav2}, its spectrum consists of a sequence of discrete simple eigenvalues, denoted 
in nondecreasing modulus order by $\lambda_n = \lambda_n(m,\theta)$, $\displaystyle{|\lambda_n|\rightarrow+\infty}$.
The associated spectral projections are of rank $1$, and E.-B. Davies showed in \cite{Dav2}
that 
for all $m\in]0,+\infty[$ and $\theta\neq0$ satisfying (\ref{argz}), for all $\alpha>0$,
there exists $N = N(m,\theta,\alpha)\geq0$ such that the instability indices $\kappa_n(m,\theta)$ 
of $\CA(m,\theta)$ satisfy
$
 \kappa_n(m,\theta)\geq n^\alpha
$
for $n\geq N$. This statement has been improved in the case $m=2$ of the harmonic oscillator (sometimes
refered as the Davies operator),
since E.-B. Davies and A. Kuijlaars showed \cite{DavKui} that $\kappa_n(2,\theta)$ grows exponentially fast as 
$n\rightarrow+\infty$, with an explicite rate $c(\theta)$ : there exists an explicit $c(\theta)$ such that
\beq\label{ResDavKui}
 \lim\limits_{n\rightarrow+\infty}\frac{1}{n}\log \kappa_n(2,\theta) = c(\theta).
\eeq
The purpose of this Note is to prove that this statement actually holds for the so-called
\emph{complex Airy operator} $\CA(1,\theta)$ 
and for the even anharmonic oscillators $\CA(2k,\theta)$, $k\geq1$.\\

\textbf{Acknowledgments}\\
I am greatly indebted to Professor Bernard Helffer for his help, advice and comments ; I am also grateful to Thierry Ramond,
Christian Gérard and André Martinez for their valuable discussions. I acknowledge the support of the ANR NOSEVOL.

\section{Non-selfadjoint anharmonic oscillators}
We first deal with the complex Airy operator $\CA(1,\theta)$ defined in (\ref{oscanh}). We show that the corresponding
instability indices $\kappa_n(1,\theta)$ grow like in
(\ref{ResDavKui}) as $n\rightarrow+\infty$.
More precisely, we get asymptotic expansions in powers of $n^{-1}$ as 
$n\rightarrow+\infty$.\\
\begin{theorem}\label{AiryGal}
Let $0<|\theta|<3\pi/4$. There exists a real sequence $(\alpha_j(\theta))_{j\geq1}$ such that the instability indices 
$\kappa_n(1,\theta)$ of $\CA(1,\theta)$ satisfy, as $n\rightarrow+\infty$,
\beq\label{equivAiD}
\exp\left(-C(\theta)(n-1/2)\right)\kappa_n(1,\theta) = \frac{K(\theta)}{\sqrt{n}}
 \left(1+\sum_{j=1}^{+\infty}\alpha_j(\theta)n^{-j}\right)+\mathcal{O}(n^{-\infty}),
\eeq
where
\[
 C(\theta) = \pi m_\theta^{3/2}|\sin\theta| ~~~ \textrm{ and } ~~~ K(\theta) = \frac{1}{2\sqrt{3|\sin\theta|}m_\theta^{1/4}},
\]
with
\[
m_\theta = \sqrt{1+\frac{\sin^2(2\theta/3)}{\sin^2\theta}-2\frac{\cos(\theta/3)\sin(2\theta/3)}{\sin\theta}} > 0.
\]
\end{theorem}
\textbf{Sketch of the proof : }
Let us first recall that all the eigenvalues of $\CA(m,\theta)$, $m\in\mathbb{N}$,
 have associated spectral projections of rank $1$, see Lemma $5$ in \cite{Dav2}. Hence, one
can easily check that \cite{AslDav}
\beq\label{kappan}
\kappa_n(m,\theta) = \frac{\|u_n\|^2}{\sca{u_n}{\bar u_n}},
\eeq
where $u_n$ denotes an eigenfunction associated with the $n$-th eigenvalue of $\CA(m,\theta)$.\\
We get rid of the singularity of the potential at $x=0$ by decomposing $\CA(1,\theta)$ into its Dirichlet and Neumann realizations
$\CA^D(1,\theta)$ and $\CA^N(1,\theta)$ on $\mathbb{R}^+$. We then compute their instability indices
\beq\label{kappaAi}
\kappa_n^{D/N}(1,\theta) = \frac{\int_{\mathbb{R}^+}|Ai(\mu_n^{D/N}+e^{i\theta/3}x)|^2dx}
{\left|\int_{\mathbb{R}^+}Ai(\mu_n^{D/N}+e^{i\theta/3}x)^2dx\right|},
\eeq
given by formula (\ref{kappan}), where $x\mapsto Ai(\mu_n^{D/N}+e^{i\theta/3}x)$ is the $n$-th eigenfunction of
$\CA^{D/N}(1,\theta)$, $\mu_n^D$ (resp. $\mu_n^N$) being the $n$-th (negative) zero (resp. critical point) of the Airy function
$Ai$ (see \cite{Alm}, \cite{Hel}).
We estimate the numerator in (\ref{kappaAi}) by using the well-known asymptotic expansion of $Ai$
at infinity in the complex plane (see \cite{AbrSteg}), and the Laplace method brings an $\exp(c_\theta|\mu_n^{D/N}|^{3/2})$ term
in $\kappa_n^{D/N}(1,\theta)$, $c_\theta>0$. The integral in the denominator of (\ref{kappaAi}), after 
deformation of the path of integration by homotopy, is equal to
\beq\label{magique}
 \int_{\mu_n^{D/N}}^{+\infty}Ai^2(x)dx = Ai'^2(\mu_n^{D/N})
\eeq
(it is indeed straightforward, using Airy equation, to check that $x\mapsto xAi^2(x)-Ai'^2(x)$ is a primitive for $Ai^2$).
Hence the expansion of $Ai'(-z)$ as $z\rightarrow+\infty$, given in \cite{AbrSteg}, provides an asymptotic expansion for
(\ref{magique}) in powers of $|\mu_n^{D/N}|^{-3/2}$. The statement follows from the behavior of 
$\mu_n^{D/N}$ as $n\rightarrow+\infty$, since we have asymptotic expansions for $(n-1/4)^{-2/3}|\mu_n^D|$
(resp. $(n-3/4)^{-2/3}|\mu_n^N|$) in powers of $(n-1/4)^{-2}$ (resp. $(n-3/4)^{-2}$).
\hfill $\square$\\


Notice that the exponential instability appears as soon as $\theta\neq0$.\\
We have a similar statement for even anharmonic oscillators :
\begin{theorem}\label{equivAnharm}
 Let $k\in\mathbb{N}^*$ and $\theta$ be such that $0<|\theta|<\frac{(k+1)\pi}{2k}$.
If $\kappa_n(2k,\theta)$ denotes the $n$-th instability index of
$
 \CA(2k,\theta) = -\frac{d^2}{dx^2}+e^{i\theta}x^{2k},
$
 then there exist $K(2k,\theta)>0$ and a real sequence $(C^j(2k,\theta))_{j\geq1}$ such that
\beq\label{dvptAnh}
e^{-c_k(\theta) n}\kappa_n(2k,\theta)=\frac{K(2k,\theta)}{\sqrt{n}}
\left(1+\sum_{j=1}^{+\infty}C^j(2k,\theta)n^{-j}\right) + \mathcal{O}(n^{-\infty})
\eeq
as $n\rightarrow+\infty$, with
\beq\label{ck}
c_k(\theta) = \frac{2(k+1)\sqrt{\pi}\Gamma\left(\frac{k+1}{2k}\right)\varphi_k(\xi_k)}{\Gamma\left(\frac{1}{2k}\right)} >0,
\eeq
where
\begin{eqnarray}
 \xi_k &  = & \left(\frac{\tan (|\theta|/(k+1))}{\sin (k|\theta|/(k+1))+\cos(k|\theta|/(k+1))
\tan (|\theta|/(k+1))}\right)^{\frac{1}{2k}},\\
\varphi_k(\xi) & = & \Im\int_0^{\xi e^{i\frac{\theta}{2(k+1)}}}(1-t^{2k})^{1/2}dt.
\end{eqnarray}
\end{theorem}

\textbf{Sketch of the proof : }
We first perform an analytic dilation and a scale change to recover the semiclassical selfadjoint anharmonic oscillator
$
 \CP_h(2k) = -h^2\frac{d^2}{dx^2}+x^{2k}-1$,
with 
\[h = h_n = |\lambda_n(2k,\theta)|^{-\frac{k+1}{2k}}.\]
 The $n$-th instability index of $\CA(2k,\theta)$ then writes
\beq\label{kappaAnh}
 \kappa_n(2k,\theta) = \frac{\int_\mathbb{R}|\psi_h(e^{i\frac{\theta}{2(k+1)}}x)|^2dx}{\int_\mathbb{R}\psi_h^2(x)dx}
\eeq
where $\psi_h$ solves $\CP_h(2k)\psi_h = 0$, $\psi_h\in L^2(\mathbb{R})$
(see (\ref{kappan}), after deformation of the integration path in the denominator). The complex WKB method 
(see \cite{Olv}, \cite{Vor}, \cite{GerGri}) and the analysis of the Stokes lines of $\CP_h(2k)$ provide an
asymptotic expansion of $\psi_h(e^{i\frac{\theta}{2(k+1)}}x)$ as $h\rightarrow0$, which enables us to determine the asymptotic
behaviour of the numerator in (\ref{kappaAnh}), using again the Laplace method.\\
On the real axis, $\psi_h$ is treated separately in its oscillatory region $[-1+\delta,1-\delta]$, $\delta>0$, and in the
neighbourhood of the turning points $\pm1$. Hence, the stationary phase method leads to an asymptotic expansion 
in powers of $h$ of the
denominator in (\ref{kappaAnh}). Finally, the statement follows from the Bohr-Sommerfeld quantization rule for $h_n$
(see \cite{GriSjo}, Exercise $12.3$) or Weyl formula \cite{HelRob}.
\hfill $\square$ \\

 In the harmonic case $k=1$ (Davies operator), the first term in (\ref{dvptAnh}) yields the Davies-Kuijlaars theorem
\cite{DavKui} :
\[
 \lim\limits_{n\rightarrow+\infty}\frac{1}{n}\log\|\Pi_n\|  =  c_1(\theta) = 
4\varphi_1\left(\frac{1}{\sqrt{2\cos(\theta/2)}}\right)
 = 2\Re f\left(\frac{e^{i\theta/4}}{\sqrt{2\cos(\theta/2)}}\right)
\]
where
$
 f(z) = \log(z+\sqrt{z^2-1})-z\sqrt{z^2-1}$.

\section{Eigenfunctions and semigroups}
The following theorem has been proved in \cite{Alm} in the
case of complex Airy operator $\CA(1,\theta)$, and in \cite{Dav2} in the harmonic case ($k=1$),
 as well as for $\CA(2k,\theta)$, 
$k\geq2$, $|\theta|<\frac{\pi}{2}$. The proof actually extends to any
operator $\CA(2k,\theta)$ with $|\theta|<\frac{(k+1)\pi}{2k}$ : 

\begin{theorem}\label{total}
For any $m=1,2k$, $k\geq1$, and any $\theta$ satisfying (\ref{argz}),
 the eigenfunctions of $\CA(m,\theta)$
 form a complete set of the space $L^2(\mathbb{R})$.
\end{theorem}

Notice however that the eigenfunctions of $\CA(1,\theta)$ and $\CA(2k,\theta)$, $k\geq1$,
can not form a Riesz basis because of the growth of the instability indices as $n\rightarrow+\infty$.\\
Theorem \ref{total} and the previous estimates enable us 
to study the convergence of the operator series defining the semigroup
$e^{-t\CA(m,\theta)}$ associated with $\CA(m,\theta)$
when decomposed along the projections $\Pi_n$.\\
The following statement extends the result of \cite{DavKui} in the harmonic case.

\begin{corollary}\label{seriesemigpes}
Let $|\theta|\leq\pi/2$ and $e^{-t\CA(m,\theta)}$ be the semigroup generated by $\CA(m,\theta)$,
$\displaystyle{\lambda_n = \lambda_n(m,\theta)}$ the eigenvalues of $\CA(m,\theta)$, and $\Pi_n = \Pi_n(m,\theta)$ the
associated spectral projections.\\
Let $\displaystyle{T(\theta) = c_1(\theta)/\cos(\theta/2)}$, where $c_1(\theta)$ is the constant in (\ref{ck}).
The series
\[
 \Sigma_{m,\theta}(t) = \sum_{n=1}^{+\infty}e^{-t\lambda_n(m,\theta)}\Pi_n(m,\theta)
\]
is not normally convergent in cases $m=1$ for any $t>0$, and $m=2$ for $t<T(\theta)$ ;
in cases $m=2$ for $t>T(\theta)$, and $m=2k$ for any $t>0$, $k\geq2$,
 the series converges normally towards $e^{-t\CA(m,\theta)}$ and,
 for $N$ sufficiently large and for some constants $C_1=C_1(k,\theta)$ and $C_2=C_2(\theta)$, the following
estimate holds :
\beq\label{restesemigpe}
 \|e^{-t\CA(m,\theta)}(I-\Pi_{<N})\| \leq 
\left\{\begin{array}{cc}
\frac{C_1}{\sqrt{N}}e^{c_k(\theta)n}\exp(-t\Re\lambda_N), & k\geq2 \\
\frac{C_2}{\sqrt{N}}\exp(-2\cos(\theta/2)(t-T(\theta))N), & k=1,~t>T
\end{array}\right.
\eeq
where $\Pi_{<N} =\Pi_1+\dots+\Pi_{N-1}$ denote the projection on the first $N-1$ eigenspaces.
\end{corollary}


\end{document}